 \documentclass[final,5p,times,twocolumn,authoryear]{elsarticle}

\usepackage{amssymb}
\usepackage{lipsum}
\usepackage{graphicx} %插入图片
\graphicspath{{Figures/}} %文章所用图片在当前目录下的Figures目录
\usepackage{booktabs} %画三线表
\usepackage{amsmath,amsthm,amssymb,amsfonts} %数学公式包

\usepackage{natbib}
\setcitestyle{numbers,square}
\newcommand{\upcite}[1]{{\setcitestyle{square,super}\cite{#1}}}
\usepackage{hyperref} %生成引用链接,注：该宏包可能与其他宏包冲突,故放在所有引用的宏包之后
\usepackage{cleveref} %实现图片和表格、公式的引用,cleveref包必须放在hyperref包之后

\usepackage{subcaption}

\begin{document}

\begin{frontmatter}

%% Title, authors and addresses

%% use the tnoteref command within \title for footnotes;
%% use the tnotetext command for theassociated footnote;
%% use the fnref command within \author or \affiliation for footnotes;
%% use the fntext command for theassociated footnote;
%% use the corref command within \author for corresponding author footnotes;
%% use the cortext command for theassociated footnote;
%% use the ead command for the email address,
%% and the form \ead[url] for the home page:
%% \title{Title\tnoteref{label1}}
%% \tnotetext[label1]{}
%% \author{Name\corref{cor1}\fnref{label2}}
%% \ead{email address}
%% \ead[url]{home page}
%% \fntext[label2]{}
%% \cortext[cor1]{}
%% \affiliation{organization={},
%%            addressline={}, 
%%            city={},
%%            postcode={}, 
%%            state={},
%%            country={}}
%% \fntext[label3]{}

\title{Physical Information Neural Networks for Solving High-index Differential-algebraic Equation Systems Based on Radau Methods}

%% use optional labels to link authors explicitly to addresses:
%% \author[label1,label2]{}
%% \affiliation[label1]{organization={},
%%             addressline={},
%%             city={},
%%             postcode={},
%%             state={},
%%             country={}}
%%
%% \affiliation[label2]{organization={},
%%             addressline={},
%%             city={},
%%             postcode={},
%%             state={},
%%             country={}}

\author[first]{Jiasheng Chen}
% % \affiliation[first]{Institute of Computing Science and Technology, Guangzhou University,  Guangzhou, China}
\author[first,second,third,fourth]{Juan Tang\corref{mycorrespondingauthor}}
\cortext[mycorrespondingauthor]{Corresponding author}
\ead{tangjn16@gzhu.edu.cn}
\author[third,fourth]{Ming Yan}
\author[second]{Shuai Lai}
\author[first]{Kun Liang}
\author[fifth]{Jianguang Lu}
\author[sixth]{Wenqiang Yang}

\address[first]{Institute of Computing Science and Technology, Guangzhou University, Guangzhou, China
}
\address[second]{School of Computer Science and Cyber Engineering, Guangzhou University, Guangzhou, China
}
\address[third]{Institute of High Performance Computing (IHPC), Agency for Science, Technology and Research (A*STAR), Singapore
}
\address[fourth]{Centre for Frontier AI Research (CFAR), Agency for Science, Technology and Research (A*STAR), Singapore
}
\address[fifth]{State Key Laboratory of Public Big Data, Guizhou University, Guiyang, China
}
\address[sixth]{Chongqing Institute of Green and Intelligent Technology, Chinese Academy of Sciences, Chongqing, China
}

\begin{abstract}
%% Text of abstract
As is well known, differential algebraic equations (DAEs), which are able to describe dynamic changes and underlying constraints, have been widely applied in engineering fields such as fluid dynamics, multi-body dynamics, mechanical systems and control theory. In practical physical modeling within these domains, the systems often generate high-index DAEs. Classical implicit numerical methods typically result in varying order reduction of numerical accuracy when solving high-index systems.~Recently, the physics-informed neural network (PINN) has gained attention for solving DAE systems. However, it faces challenges like the inability to directly solve high-index systems, lower predictive accuracy, and weaker generalization capabilities. In this paper, we propose a PINN computational framework, combined Radau IIA numerical method with a neural network structure via the attention mechanisms, to directly solve high-index DAEs. Furthermore, we employ a domain decomposition strategy to enhance solution accuracy. We conduct numerical experiments with two classical high-index systems as illustrative examples, investigating how different orders of the Radau IIA method affect the accuracy of neural network solutions. The experimental results demonstrate that the PINN based on a 5th-order Radau IIA method achieves the highest level of system accuracy. Specifically, the absolute errors for all differential variables remains as low as $10^{-6}$, and the absolute errors for algebraic variables is maintained at $10^{-5}$, surpassing the results found in existing literature. Therefore, our method exhibits excellent computational accuracy and strong generalization capabilities, providing a feasible approach for the high-precision solution of larger-scale DAEs with higher indices or challenging high-dimensional partial differential algebraic equation systems.
\end{abstract}

%%Graphical abstract
%\begin{graphicalabstract}
%\includegraphics{grabs}
%\end{graphicalabstract}

%%Research highlights
%\begin{highlights}
%\item Research highlight 1
%\item Research highlight 2
%\end{highlights}

\begin{keyword}
%% keywords here, in the form: keyword \sep keyword, up to a maximum of 6 keywords
Differential algebraic equation \sep Radau IIA method \sep Physics-informed neural network \sep Domain decomposition

%% PACS codes here, in the form: \PACS code \sep code

%% MSC codes here, in the form: \MSC code \sep code
%% or \MSC[2008] code \sep code (2000 is the default)

\end{keyword}

\end{frontmatter}

%\tableofcontents

%% \linenumbers

%% main text

\section{Introduction}
\label{introduction}

The concept of differential algebraic equations (DAEs) was formally proposed by Gear in the study of network analysis and continuous system simulation problems.\upcite{Gear1971}~Petzold made it explicit through his study of numerical methods that DAEs are not ordinary differential equations (ODEs).\upcite{Petzold1982}~DAE systems are composed of coupled ODE systems and algebraic equation systems with physical significance. These systems encompass both differential and algebraic variables, and their system form is more generalized compared to traditional ODE systems. DAEs have gained significant attention since their inception, as they can accurately describe systems that some ODEs cannot represent. They have found extensive applications in various fields, including fluid dynamics, multi-body dynamics, electronic circuits, mechanical systems, control theory, and chemical engineering.
\par
In different developmental periods and research fields, DAEs are also known as singular systems, general systems, descriptor systems, or constrained systems, among other names. They often exhibit various structural forms, such as linear DAEs, nonlinear DAEs, semi-explicit DAEs, implicit DAEs, and Hessenberg-type DAEs. Fortunately, in practical physical modeling, most of the system models obtained are either low-index DAEs or high-index ($>= $2) Hessenberg-type DAEs.\upcite{Ascher1998}~ The index of DAEs measures the 'distance' between DAEs and ODEs. Generally, a higher index implies greater difficulty in transforming DAEs into ODEs or in directly solving DAEs using ODE numerical methods. Traditional numerical methods for solving DAE systems include implicit Runge-Kutta methods\upcite{Ascher1991}, BDF methods\upcite{Cash2000}, pseudospectral methods \upcite{Saravi2010}, adomian decomposition method \upcite{Hosseini2006}, exponential integrators \upcite{Newman2003}, generalized-$\alpha$ methods \upcite{Ding2013}, and Lie group methods\upcite{Lu2016,Liu2017,Tang2023}
. It's worth noting that these direct numerical methods can solve DAEs with an index of 1. However, for high-index DAE systems, these methods are only applicable to a certain class of DAEs and may result in varying order reduction of numerical accuracy.
\par
With the rapid advancement of neural network technology and hardware resources, neural networks are demonstrating increasingly powerful capabilities. Compared to traditional numerical computing methods, neural networks offer several advantages, including strong generalization, fault tolerance, and the ability for parallel computation.~In 1998, Lagaris et al. \upcite{Lagaris1998} approximated solutions to ODEs or PDEs problems by constructing parameterized trial functions. These trial functions consist of two parts: one part satisfies initial conditions or boundary conditions which does not contain trainable parameters, while the other part is a simple feed forward neural network with trainable parameters.~In 2019, Raissi et al.\upcite{Raissi2019} introduced an important technique known as Physics-informed Neural Network (PINN) for the numerical approximation of partial algebraic equations (PDEs) problems. The PINN loss function includes not only initial or boundary conditions that reflect physical properties but also a residual term at selected points in the time-space domain where the PDEs hold. It's worth noting that PINN is a data-driven approach that doesn't require prior knowledge of the analytical form of the solution; instead, it learns the solution from data. Various variants of PINN have been proposed based on different collocation methods, such as variational hp-VPINN \upcite{Kharazmi2021} and conservative PINN (CPINN) \upcite{Jagtap2020}. Additionally,~PINN has been widely applied to solve problems in various fields, including fluid dynamics \upcite{Mao2020,Jin2021}, seismic wave prediction \upcite{Borate2023}, and optical problems \upcite{Chen2022}.
 \par
 In recent years, many researchers have attempted to construct neural network models from different perspectives to solve various types of DAEs systems influenced by these methods. For Hessenberg-DAEs with control variables and an index of $3$, Kozlov and Tiumentsev \upcite{Kozlov2015} achieved the implementation of BDFs method using a semi-empirical neural network model. Zhao Yang et al.\upcite{Yang2019} constructed a single-layer feed-forward neural network (FFNN) to solve Hessenberg-type DAEs systems. They augmented the loss function in their special Euler-Lagrange equation system with penalty terms for algebraic equations to avoid drifting in the results. Experimental results in their paper showed that the FFNN method with Sigmoid activation function provided approximate analytical solutions close to the numerical solutions of corresponding Runge-Kutta methods, but they didn't provide further details about the method's accuracy.~For linear DAEs systems, Hongliang Liu et al.\upcite{Liu2021} selected Jacobi polynomials as activation functions and constructed a single-hidden-layer feed-forward neural network (JNN). They determined the network parameters using the classical ELM algorithm. Through experimental comparisons with other approximation methods such as Pad\'{e} approximation, ADM method, and Adams methods, they illustrated the feasibility and superiority of the JNN method. It's worth noting that the examples in the paper involve DAEs with an index of 1 or linear DAEs that have been reduced to index 1. For DAEs systems with an index of 1, Moya et al. \upcite{Moya2023} proposed a neural network architecture called DAE-PINN based on the PINN method for solving DAEs systems. This neural network model is a discrete-time model based on the implicit Runge-Kutta method, which can directly address most index-1 differential-algebraic equation problems. However, it cannot solve high-index DAEs problems and suffers from low accuracy issues. To address the high-accuracy computation challenges in high-index DAEs systems, we have combined the Radau IIA numerical method with a neural network structure based on attention mechanisms. We have proposed a PINN computational framework based on the Radau method. Furthermore, we have improved the efficiency and accuracy of the solution by applying a strategy of domain decomposition.
\par
In section $2$, we briefly introduce the fundamental concepts of DAEs systems, the Radau IIA numerical method, and the neural network structure based on attention mechanisms. Building upon this foundation, we provide a detailed construction of the PINN computing framework based on the Radau IIA method. Additionally, we employ a time domain decomposition strategy for neural network.~Section $3$ use the neural network designed in this paper to solve two high-index DAEs systems, and we analyze the solving accuracy of this neural network.~Finally, we discuss and summarize the advantages, challenges, and potential avenues for improvement in the Radau-PINN architecture.

\section{Scientific Machine Learning Methods}
%%\label{}
%\lipsum[1]
This section first sequentially introduces the basic concepts of DAEs and the classical Radau IIA numerical method. Then, we introduce a neural network structure based on attention mechanisms. Building upon this, we construct a PINN based on the Radau IIA method. Finally, we enhance the efficiency and accuracy of neural network solutions for DAEs systems by utilizing the concept of domain decomposition. 
\subsection{Radau IIA Method for DAE Systems}
This article first provides a brief introduction to DAEs with an index of $2$, with the specific form as follows:
\begin{eqnarray}
\left\{\begin{aligned}y'(t) = & f(t, y(t), z(t)),\\
0 = & g(t,y(t)),
\end{aligned}
\right.
\end{eqnarray}
where~$y(t)\in {\rm{R}^n}$ is the differential function variable,~$z(t) \in {\rm{R}^m}$ is the algebraic function variable,~$t \in \left[ {{t_0},T} \right]$,~${t_0}$ is the initial time point, and~$y_0=y(t_0)$ is the initial value. Both~$f(t,y,z) \in {\rm{R}^n}$ and $g(t,y) \in {\rm{R}^m}$ are sufficiently smooth, and the Jacobian matrix $g_y f_z$ is non-singular.

The Radau IIA method is a class of implicit Runge-Kutta methods, typically defined in the following general form:
\begin{eqnarray}\label{eq:Radau}
	{\xi _i} = {y_n} + h\mathop \sum \limits_{j = 1}^v {a_{i,j}}f\left( {{\xi _j},{\zeta _j}} \right),  \\
	g\left( {{\xi _i},{\zeta _i}} \right) = 0, \\
	{y_{n + 1}} = {y_n} + h\mathop \sum \limits_{j = 1}^v {b_j}f\left( {{\xi _j},{\zeta _j}} \right),  \\
	g\left( {{y_{n + 1}},{z_{n + 1}}} \right) = 0,
\end{eqnarray}
where~${\xi _i} = y\left( {{t_n} + {c_i}h} \right)$,~${\zeta _i} = z\left( {{t_n} + {c_i}h} \right)$,~$h$ is the step size,~$n$ is the current step number,~$\left\{ {{a_{ij}},{b_j},{c_i}} \right\}$ are parameters, and $c_i = \sum\limits_{j = 1}^v {{a_{ij}}}$,~$i,j = 1, \cdots,v$.

\begin{table}[h]
    \centering
$\begin{array}{c|ccccc}
	c_1 & a_{11} & a_{12} & a_{13} & \cdots & a_{1v} \\
	c_2 & a_{21} & a_{22} & a_{23} & \cdots & a_{2v} \\
	c_3 & a_{31} & a_{32} & a_{33} & \cdots & a_{3v} \\
	\vdots & \vdots & \vdots & \vdots & \ddots & \vdots \\
	c_v & a_{v1} & a_{v2} & a_{v3} & \cdots & a_{vv} \\
	\hline
	& b_1 & b_2 & b_3 & \cdots & b_v \\
\end{array}$
    \caption{The parameter table of the $v$-stage implicit Runge-Kutta methods.
}
\label{tab:RK}
\end{table}

In table \ref{tab:RK}, different sets of parameters lead to different implicit Runge-Kutta methods, such as commonly used Gauss method, Radau method, and Lobatto method. These parameters are determined using Gauss polynomials, Radau polynomials, and Lobatto polynomials, respectively. Among them, the Radau IIA method is a high-precision numerical method with excellent numerical stability. Therefore, in this paper, the Radau IIA method is chosen, and the parameters need to satisfy the following conditions:
\begin{align}
B(2v-1):& \sum_{i=1}^v b_i c_i^{k-1}= \frac{1}{k},~k=1,\cdots,2v-1 \\
C(v):& \sum_{j=1}^v a_{ij} c_j^{k-1}= \frac{c_i^k}{k},~k=1,\cdots,v  \\
D(v-1):& \sum_{i=1}^v b_i c_i^{k-1} a_{ij}= \frac{b_j}{k}(1-c_j^k),~k=1,\cdots,v-1
\end{align}
and $c_v=1$,~$b_j=a_{vj}$,~$i,j=1,2,\dots,v$.

\subsection{Neural Network Structure Based on Attention Mechanism}
\begin{figure*}[!h]
	\centering 
	\includegraphics[width=0.85\linewidth]{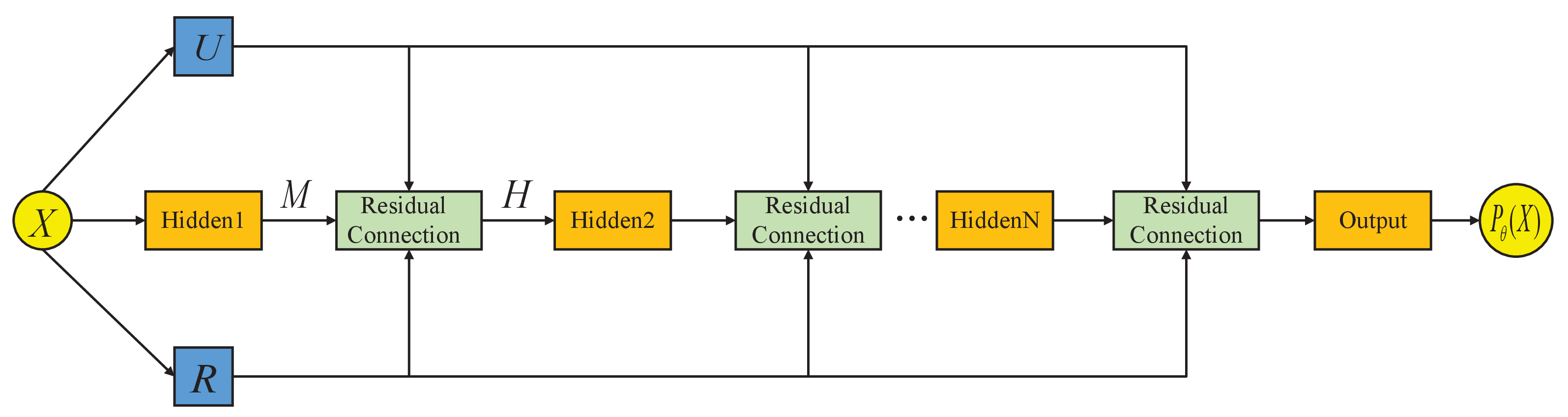}	
	\caption{Improved attention neural network structure.} 
	\label{fig:Atten}
\end{figure*}

%%\label{}
%\lipsum[2]
Building upon the DAE-PINN structure, we employ adaptive activation functions (\ref{eq:sig}) to train a neural network structure based on an attention mechanism. The specifics are as follows:

The improved neural network model architecture based on attention mechanisms is primarily constructed using two Transformer networks, denoted as U and R, to build two stacked layer networks, as illustrated in Figure \ref{fig:Atten}. Both neural networks map the input variable $X$ (differential function variable $y$) to a high-dimensional feature space. Subsequently, each hidden layer forms new residual connections using element-wise multiplication operations, as expressed below:
\begin{eqnarray}\label{eq:sig}
U = & \phi (X{W^1} + {b^1}), \\
R =& \phi (X{W^2} + {b^2}), \\
{H^{(1)}} =& \phi (\eta  \cdot l \cdot X{W^{o,1}} + {b^{o,1}}),  \\
{M^{({\rm{k}})}} =& \phi ({H^{{\rm{k}}}}{W^{o,k}} + {b^{o,k}}), \\
{H^{(k + 1)}} =& (1 - {M^{(k)}}) \odot U + {M^{(k)}} \odot R, \\
{P_\theta }(X) =& {H^{{\rm{d + 1}}}}W + b,
\end{eqnarray}
where~$X$ represents the input vector of the neural network, ~${W^{o,k}}$ is the collection of weights for the $o$-th neuron in the $k$-th layer, ~${b^{o,k}}$ denotes the set of biases for the $o$-th neuron in the $k$-th layer, ~$\phi$ is the activation function, ~$\odot$ represents element-wise multiplication, $d$ indicates the number of hidden layers (the depth of the neural network), ~${P_\theta }(X)$ is the final output vector of the neural network, ~$\eta$ is a predetermined hyper-parameter that ensures the slope is greater than 1, and $l$ is a parameter that can modify the slope of the activation function.

\begin{figure*}[!h]
	\centering 
	\includegraphics[width=1\linewidth]{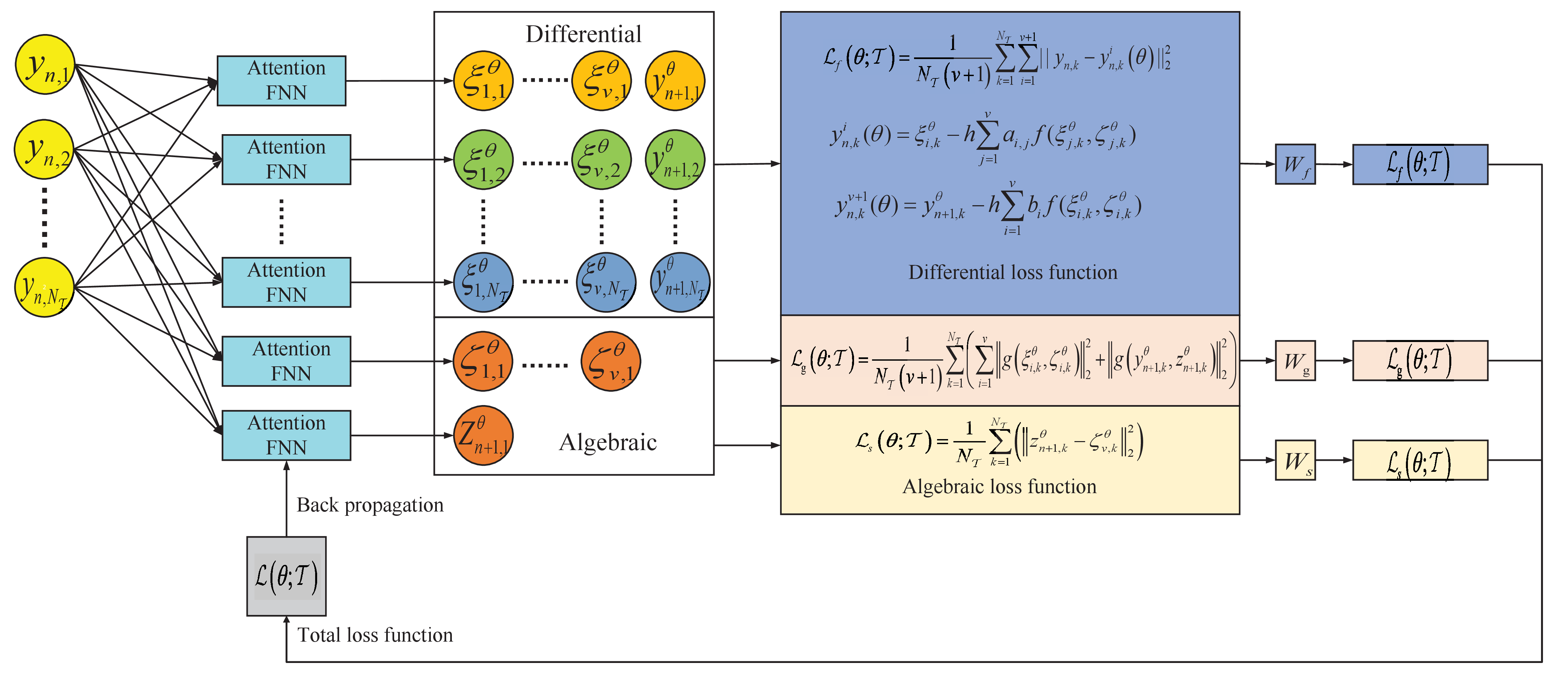}	
	\caption{The schematic diagram of PINN based on Radau IIA method.} 
	\label{fig:Radau-PINN}
\end{figure*}

\begin{figure*}[t]
	\centering 
	\includegraphics[width=0.70\linewidth]{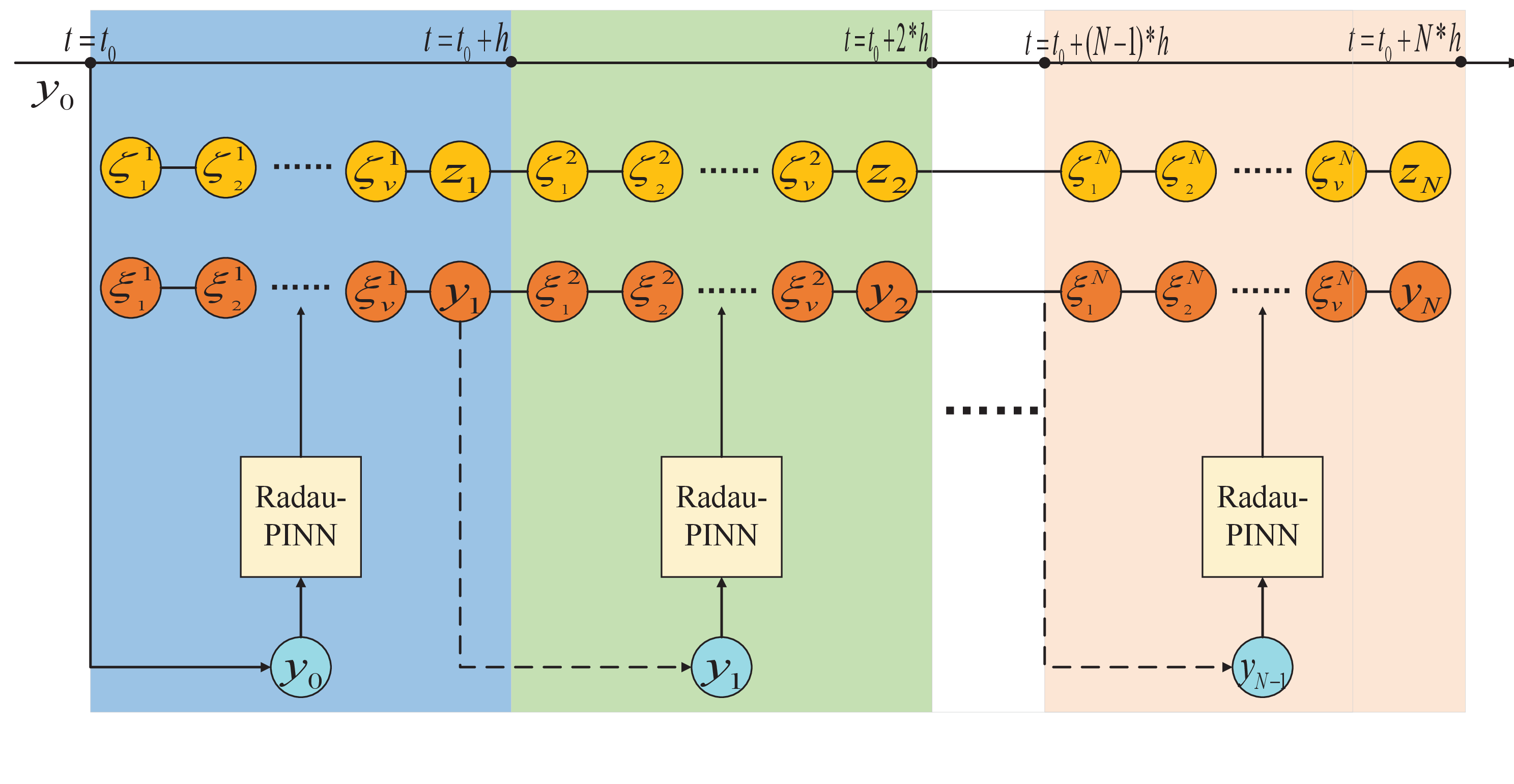}	
	\caption{The time domain decomposition of neural networks.} 
	\label{fig:Time}
\end{figure*}

\subsection{PINN Based on Radau IIA Method}
%\lipsum[3]
In this section, we use the discrete-time model of PINN as the foundation, incorporating a neural network structure based on attention mechanisms. We have constructed a PINN architecture based on the Radau IIA method, as illustrated in Figure \ref{fig:Radau-PINN}.

Firstly, construct a neural network with multiple inputs and multiple outputs, where the inputs consist of the collection of differential variables ${y_{n}}$, and outputs
\begin{eqnarray}
\xi _1^\theta ,\xi _2^\theta ,,......,\xi _v^\theta ,{y_{n + 1}}^\theta; \\
\zeta _1^\theta ,\zeta _2^\theta ,,......,\zeta _v^\theta ,{z_{n + 1}}^\theta. 
\end{eqnarray}

The first $v$ values of $\xi _{i}^\theta $ represent intermediate differential variables, and the first $v$ values of $\zeta _{i}^\theta $ represent intermediate algebraic variables, where ~$i = 1,2,......,v$.

Secondly, based on the structure of the DAEs system with an index of 2 and the characteristics of the Radau IIA method, we further design neural network structures based on the attention mechanism for both the differential variable part and the algebraic variable part of the system. There are two specific design approaches: one assigns a single neural network to all differential variables and two neural networks to the algebraic variables, and the other assigns a separate neural network to each differential variable while keeping the algebraic variable part unchanged. In the case of the algebraic variable part, one of the neural networks is used to predict the first v values, while the other neural network is used to predict the v+1-th value. Theoretically, the second approach (as shown in Figure \ref{fig:Radau-PINN}) constructs a neural network for each individual differential or algebraic variable, thereby improving the overall model's accuracy and generalization. Through further testing, the second approach's training results are more precise than those of the first approach, consistent with the expected results. As a result, all subsequent experiments in this paper are implemented based on the second approach.

Thirdly, based on the designed network structure, this paper constructs the loss function as follows:
\begin{eqnarray}
{\cal L}\left( {\theta ;{\cal T}} \right) = {W_f}{{\cal L}_f}\left( {\theta ;{\cal T}} \right) + {W_{\rm{g}}}{{\cal L}_{\rm{g}}}\left( {\theta ;{\cal T}} \right) + {W_s}{{\cal L}_s}\left( {\theta ;{\cal T}} \right),
\end{eqnarray}

\begin{itemize}
    \item[$\bullet$] ${{\cal L}_f}\left( {\theta ;{\cal T}} \right)$ is the loss related to the differential network and is expressed as follows:
\begin{eqnarray}
\frac{1}{{{N_{\cal T}}\left( {v + 1} \right)}}\mathop \sum \limits_{k = 1}^{{N_{\cal T}}} \mathop \sum \limits_{i = 1}^{v + 1} ||{y_{n,k}} - y_{n,k}^i\left( \theta  \right)||_2^2,\\
y_{n,k}^i(\theta ) = \xi _{i,k}^\theta  - h\sum\limits_{j = 1}^v {{a_{i,j}}f(\xi _{j,k}^\theta ,\zeta _{j,k}^\theta )}, \\
\nonumber i = 1,\cdots,v, \\
y_{n,k}^{v + 1}(\theta ) = y_{n + 1,k}^\theta  - h\sum\limits_{i = 1}^v {{b_i}} f(\xi _{i,k}^\theta ,\zeta _{i,k}^\theta );
\end{eqnarray} 

\item[$\bullet$]
${{\cal L}_{\rm{g}}}\left( {\theta ;{\cal T}} \right)$ is the loss associated with the algebraic network and can be expressed as follows:
\begin{equation} 
\frac{1}{{{N_{\cal T}}\left( {v + 1} \right)}}\mathop \sum \limits_{k = 1}^{{N_{\cal T}}} \left( {\mathop \sum \limits_{i = 1}^v \left\| {g\left( {\xi _{i,k}^\theta ,\zeta _{i,k}^\theta } \right)} \right\|_2^2 + \left\| {g\left( {y_{n + 1,k}^\theta ,z_{n + 1,k}^\theta } \right)} \right\|_2^2} \right);
\end{equation}

\item[$\bullet$]
${{\cal L}_s}\left( {\theta ;{\cal T}} \right)$ is the loss related to the last value of the controlled algebraic variable and can be expressed as follows:
\begin{equation}
\frac{1}{{{N_{\cal T}}}}\mathop \sum \limits_{k = 1}^{{N_{\cal T}}} \left( {\left\| {z_{n + 1,k}^\theta  - \zeta _{v,k}^\theta } \right\|_2^2} \right),
\end{equation}
\end{itemize}
where~${W_f}$ represents the loss weight for the differential neural network, ~${W_g}$ is the weight for the algebraic neural network, and ~${W_s}$ signifies the weight for the control of algebraic variable prediction neural network. The parameters ~$a_{i,j}$ and ~$b_i$ are specific to the Radau IIA method. ~${\cal T}$ is the total number of samples, ~${N_{\cal T}}$ is the number of training samples in the current batch, and ~$\theta$ denotes the neural network parameters. Here, ~$f$ represents the differential network, and ~$g$ represents the algebraic network. ~${y_{n,k}}$ corresponds to the sample data of the model, ~$\xi _{i,k}^\theta$ stands for the values of intermediate differential variables, and ~$\zeta _{i,k}^\theta $ signifies the values of intermediate algebraic variables. Furthermore, ~$y_{n,k}^i\left( \theta  \right)$ represents the output values of the differential neural network. The notation ~$\left\|  \cdot  \right\|_2^2$ refers to the square of the L2 norm, ~${z_{n + 1,k}^\theta }$ represents the final output of the algebraic neural network, and ~${\zeta _{n,k}^\theta }$ denotes the penultimate output of the algebraic neural network.

Finally, we use gradient descent to solve for the weights, biases, and other parameters of the PINN,
\begin{equation}
{\theta ^*} = \arg \mathop {\min }\limits_\theta  {\cal L}\left( {\theta ;{\cal T}} \right).
\end{equation}

\subsection{Time Domain Decomposition of Neural Networks}

In this section, based on an analysis of the existing limitations of the PINN architecture, we adopt a time-domain decomposition strategy using neural networks.

One limitation of the PINN model is that it exhibits relatively low accuracy in predicting solutions. This is because the inherent inaccuracies involved in solving high-dimensional non-convex optimization problems can lead to local minima, making it challenging to achieve absolute errors below ${10^{-5}}$. Another evident limitation is the high training cost .\upcite{Jagtap2020} Similarly, our proposed PINN model based on the Radau IIA method may encounter similar issues. Furthermore, the iterative format of the Radau IIA method does not fully exploit its high-precision advantages during training.

To address these issues, we propose a time-domain decomposition strategy for neural networks, as illustrated in Figure \ref{fig:Time}. With this approach, we partition the original problem into segments, which not only enhances solution accuracy but also leverages the advantages of iterative training. In other words, the predicted values from the previous time segment can serve as input values for the subsequent segment. This means that knowing the data values at the initial point $t_0$ for the first segment is sufficient to iteratively compute the solutions over the entire time domain. This approach significantly reduces the amount of required data.~Specifically, only the data at the initial point $t_0$ for a set of differential variables, denoted as $y_0$, is needed. Using the time-domain decomposition structure, we can iteratively determine the desired values within the range $[t_0, T]$. This involves information related to $(T - t_0) / h \cdot v$ data points, which reduces the need for extensive training data. On the other hand, if we can obtain the initial values for each network at every time segment, parallel training of each neural network becomes possible, significantly reducing the model training time.

\section{Numerical Experiments}

In this section, we apply PINN based on the Radau IIA method to solve two high-index DAEs systems separately and further investigate the influence of the order of the Radau IIA method on the solution results. The experiments were conducted on a Windows 10 operating system with an Intel(R) Core(TM) i7-10875H CPU @ 2.30GHz processor. We used Python 3.9 software and coded the neural network architecture using PyTorch 1.12.1, the GPU version. Additionally, this paper involves two formulas to measure the accuracy of the experiments. One is the commonly used Absolute Error (AE) formula, defined as $AE = |y_{true} - y_{pred}|$, which reflects the magnitude of the deviation between the neural network's predicted solution and the true solution. The other metric is the Mean Absolute Error (MAE) formula, defined as $MAE = \frac{1}{n} \sum_{i=1}^n |y_{true}^i - y_{pred}^i|$, used to assess the differences in accuracy among different orders of the Radau IIA method.

\subsection{Hessenberg-type DAEs System}

In this section, we explore classical Hessenberg-type DAE systems with an index of 2 that possess exact analytical solutions\upcite{Tang2023}, as follows:

\begin{eqnarray}
\left\{\begin{aligned}
y_1'(t)= & (y_3(t) y_4(t) + y_1(t) y_2(t)) y_5(t), \\
y_2'(t)= &  - y_3(t) y_4(t)^2 y_2(t)^2 y_5(t), \\
y_3'(t)  = & 2 y_3(t) y_4(t) y_1(t) y_2(t), \\
y_4'(t) = &  - y_3(t) y_4(t) y_2(t)^2, \\
0 = & y_1(t) y_4(t) - y_2(t) y_3(t),
\end{aligned}
\right.
\end{eqnarray}
where $t \in [0,1]$, and the initial values $y_0 = (1, 1, 1, 1, 1)$. The functions $y_1(t)$, $y_2(t)$, $y_3(t)$, $y_4(t)$ represent differential variables, while $y_5(t)$ is an algebraic variable. The system's exact solution expressions are $y_1(t) = e^{2t}$, $y_2(t) = e^{-t}$, $y_3(t) = e^{2t}$, $y_4(t) = e^{-t}$, and $y_5(t) = e^t$.

\begin{figure*}[!ht]
	\centering 
	\includegraphics[width=1\linewidth]{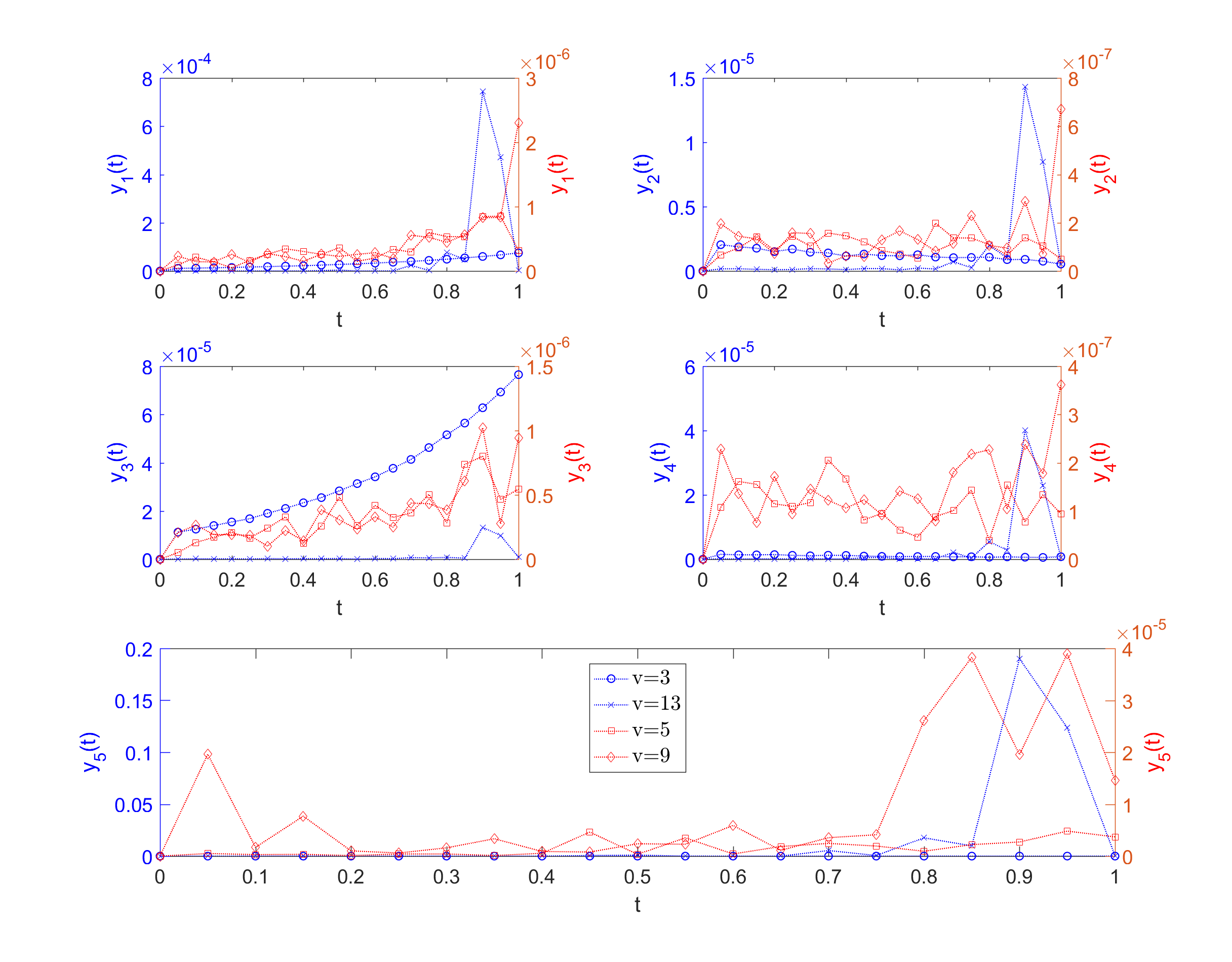}	
    \vspace{-2em}
	\caption{ The mean absolute errors for solving Hessenberg-DAEs systems using PINN based on Radau IIA of order $v=3,5,9,13$.~The blue curves represent the mean absolute errors $v=3,13$ on the left Y-axis, while the red curves correspond to $v=5,9$ on the right Y-axis.} 
 \label{fig:ex1RadauOrds2357}
\end{figure*}
\begin{figure*}[!h]
	\centering 
	\includegraphics[width=1\linewidth]{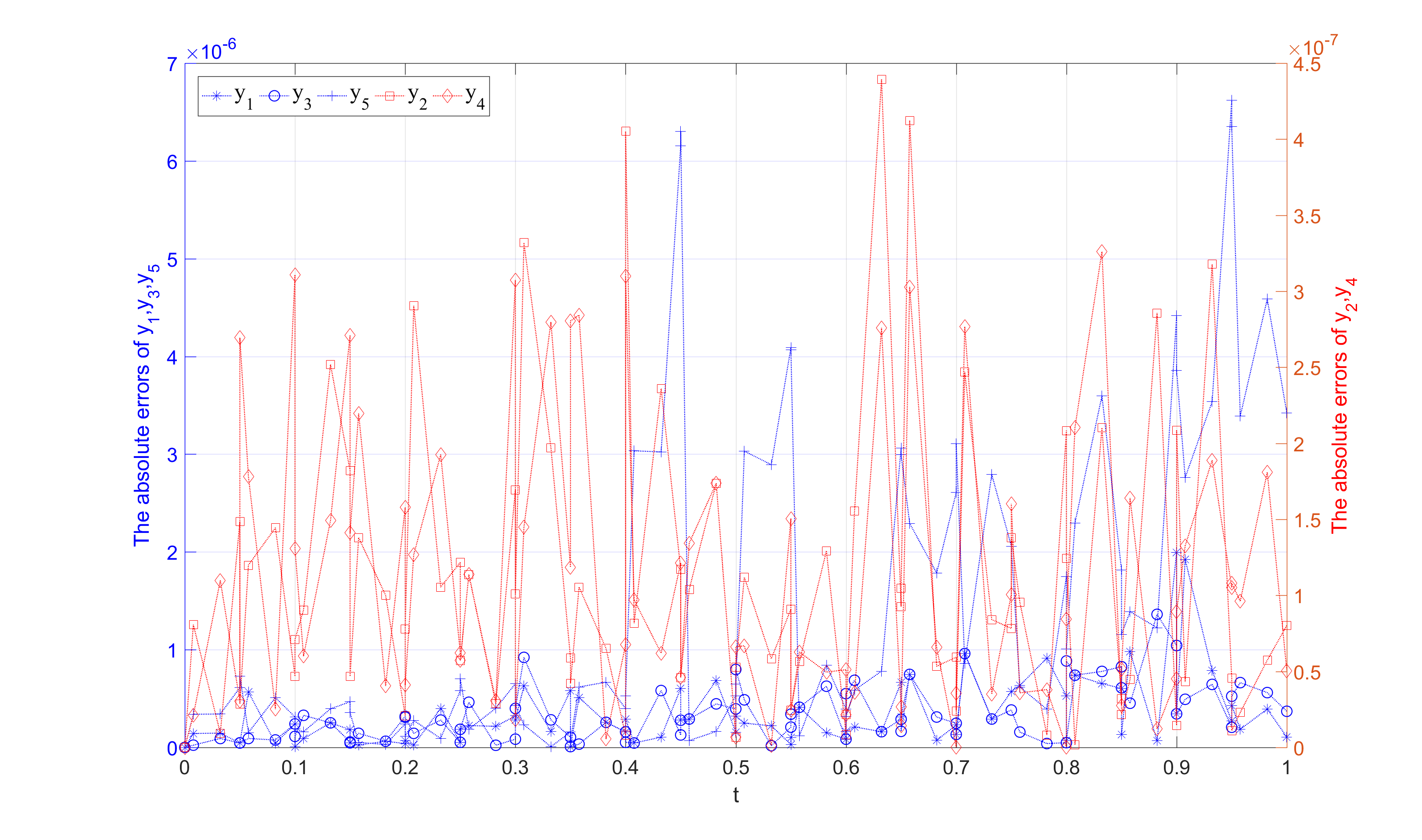}	
    \vspace{-3em}
	\caption{ The absolute errors of Hessenberg-DAEs system solved by PINN based on 5th-order Radau IIA.~The blue curves represent the absolute errors of $y_1$,~$y_3$ and $y_5$ on the left Y-axis, while the red curves correspond to $y_2$ and $y_4$ on the right Y-axis.} 
 \label{fig:ex1RadauOrds5}
\end{figure*}

Firstly, we consider the impact of different orders of Radau IIA methods, including 3rd, 5th, 9th, and 13th orders (corresponding to $v=2,3,5,7$), on the precision of neural network solutions. Secondly, we explore the influence of activation functions on PINN. Common activation functions for hidden layers include Sigmoid, TanH, Sin, and ReLu, among others. When solving smoothly continuous systems, ReLu is generally not chosen; instead, Sigmoid, TanH, or Sin activation functions are preferred. In the experiments, \textit{Sigmoid} resulted in better approximate solutions. Within this neural network framework, the initial values of the differential variables, namely $y_1(t)$, $y_2(t)$, $y_3(t)$, and $y_4(t)$ for each time segment, are used as a dataset for training. The step size $h$ is $0.05$, which means that each time interval has a length of 0.05. Each network model in every time segment comprises 5 hidden layers, with each hidden layer containing 100 neurons. Sigmoid is used as the activation function, and the Adam optimizer is applied for 100,000 iterations. The experimental results within the time interval of 0 to 1 are presented in Figures \ref{fig:ex1RadauOrds2357}.

From Figure \ref{fig:ex1RadauOrds2357}, it is evident that the accuracy of the mean absolute errors for the 3rd and 13th-order Radau IIA methods corresponds to the blue Y-axis, while the accuracy of the average absolute errors for the 5th and 9th-order Radau methods corresponds to the red Y-axis. For all the differential function variables, the 3rd and 13th-order Radau IIA methods exhibit significantly higher average absolute errors compared to the 5th and 9th-order methods. For the algebraic variable $y_5$, the 13th-order Radau IIA method has notably higher average absolute errors than the 3rd, 5th, and 9th-order methods. 

Additionally, we further observe that for all differential function variables from red Y-axis, the 9th-order Radau IIA method's overall trend in average absolute errors is significantly higher than the 5th-order method. For the algebraic variable $y_5$, the 9th-order Radau IIA method exhibits notably higher average absolute errors than the 5th-order method. In other words, the 5th-order Radau IIA-based PINN achieves the highest precision in terms of average absolute errors.

The absolute error results obtained using the 5th-order method are shown in Figure \ref{fig:ex1RadauOrds5}. The accuracy of the absolute errors for $y_1(t)$, $y_3(t)$, and $y_5(t)$ corresponds to the blue Y-axis, while the accuracy of the absolute errors for $y_2(t)$ and $y_4(t)$ corresponds to the red Y-axis. From the figure, it is evident that the neural network's predicted values for all four differential variables have their lowest precision of absolute errors maintained at the order of $10^{-6}$, while the lowest precision of absolute errors for the algebraic variable is kept at $10^{-6}$. The experimental results suggest that the neural network's predicted solutions have reached a high level of accuracy. 

For the neural network structure designed in this paper, the predicted values of the differential variables can be used as the initial values for the next time step's network input dataset. The precision of the differential variables can affect the results of the next time step's network. In this context, the precision of the differential variables $y_1$ and $y_3$ is already at the order of $10^{-6}$, and the precision of the differential variables $y_2$ and $y_4$ is at the order of $10^{-7}$, which will not significantly affect the precision of the next time step.

\subsection{DAE System of the Pendulum Model}

\begin{figure*}[!h]
	\centering 
	\includegraphics[width=1\linewidth]{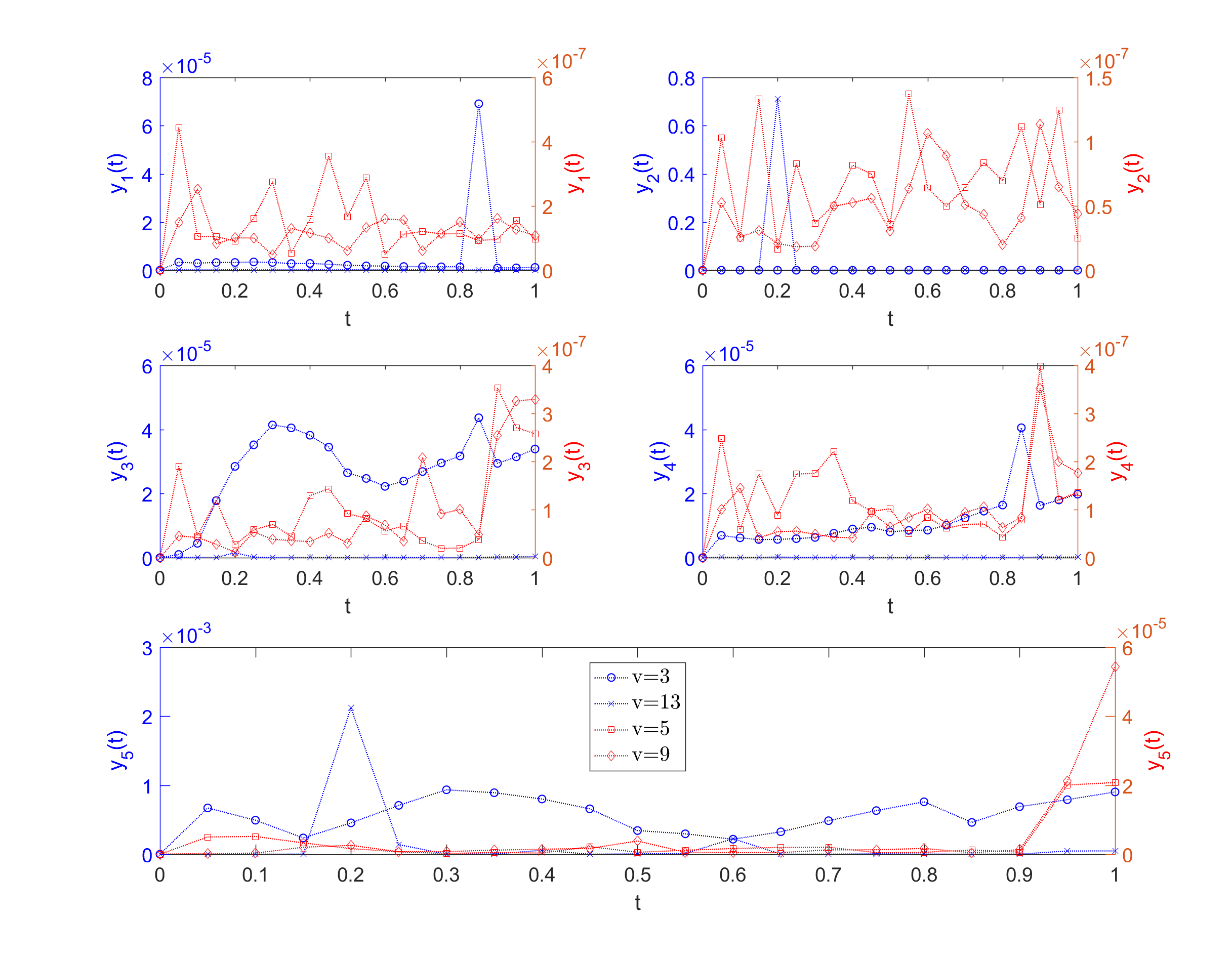}	
	\caption{ The mean absolute errors for solving single Pendulum DAEs systems using PINN based on Radau IIA of order $v=3,5,9,13$.~~The blue curves represent the mean absolute errors $v=3,13$ on the left Y-axis, while the red curves correspond to $v=5,9$ on the right Y-axis.} 
  \label{fig:ex2RadauOrds2357}
\end{figure*}

In this section, we study the classical pendulum DAEs system with an index of 2, as follows:

\begin{eqnarray}
\left\{\begin{aligned}
y_1'(t)  = & y_3(t), \\
y_2'(t)  = & y_4(t), \\
y_3'(t)  = & - y_1(t)y_5(t), \\
my_4'(t) = &  - y_2(t)y_5(t)  - \lambda,     \\
0        = & y_1(t)y_3(t) + y_2(t)y_4(t),
\end{aligned}
\right.
\end{eqnarray} 
where $t \in [0,1]$, and the parameters $m$ and $\lambda$ are variable parameters, both set to 1 in the experiments of this section. The initial values are $y_0 = (1, 0, 0, 1, 1)$. In this context, $y_1(t)$, $y_2(t)$, $y_3(t)$, and $y_4(t)$ are differential function variables, while $y_5(t)$ is an algebraic function variable. This DAEs system does not have an exact analytical expression. In this paper, we directly solve the reduced inner ODEs of this system using high-precision ODE solvers from the Python scientific computing library \textit{Scipy} and compare the obtained approximate solution with the predicted values from the neural network. Similarly, we consider the impact of different orders (3, 5, 9, 13, corresponding to $v=2,3,5,7$) in the Radau IIA methods on the accuracy of the neural network's solutions. Secondly, we explore the effect of activation functions on PINN. In this experiment, the \textit{Sin} activation function provides a better approximation. To maintain consistency in the numerical experiments, other network structural information is consistent with the experiments in the previous section. The results obtained are shown in Figures \ref{fig:ex2RadauOrds2357}.

\begin{figure*}[!h]
	\centering 
	\includegraphics[width=1\linewidth]{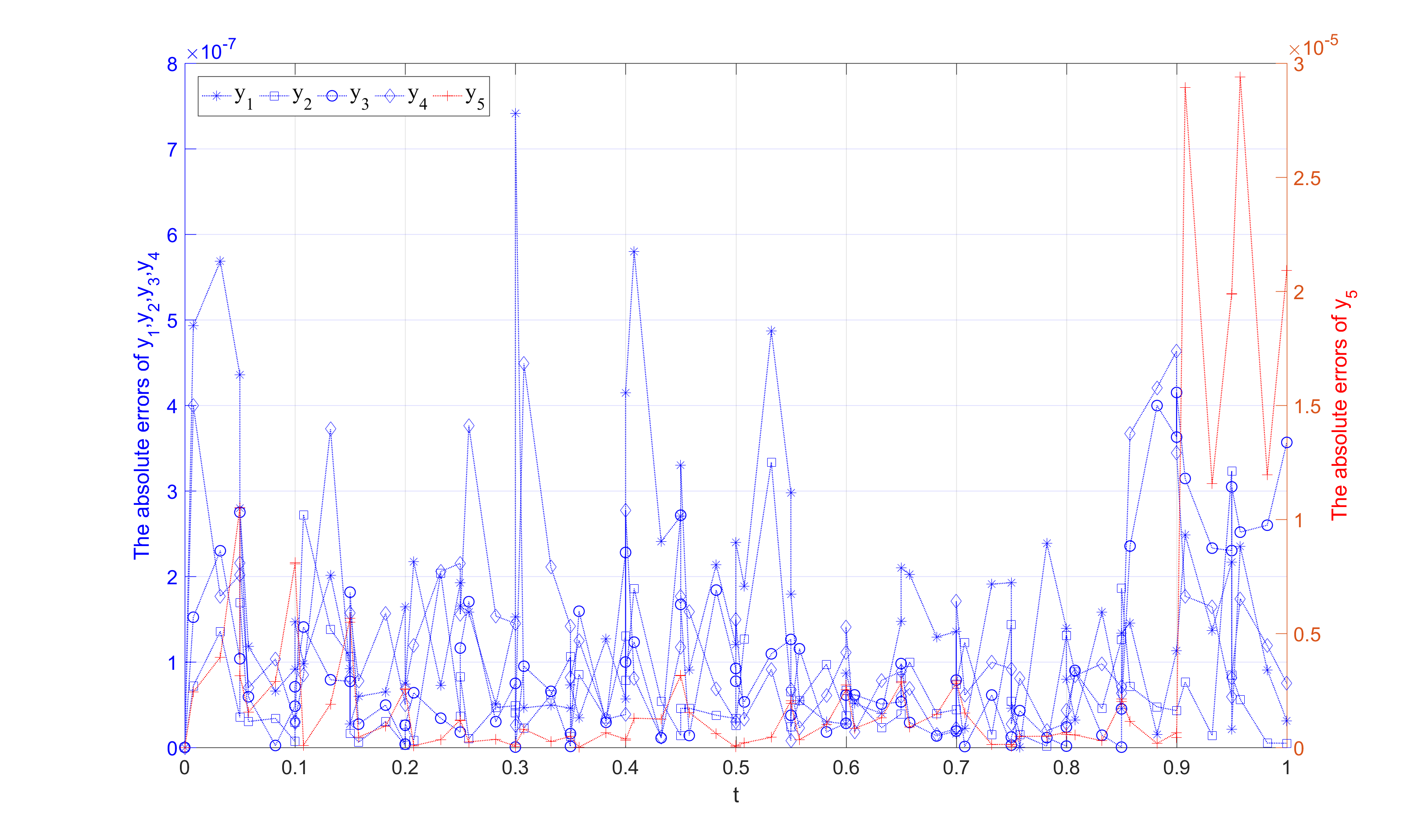}	
    \vspace{-3em}
	\caption{ The absolute errors of single Pendulum DAEs system solved by PINN based on 5th-order Radau IIA.~The blue curves represent the absolute errors of $y_1$,~$y_2$,~$y_3$ and $y_4$ on the left Y-axis, while the red curve corresponds to $y_5$ on the right Y-axis.} 
\label{fig:ex2RadauOrds5}
\end{figure*}

From Figure \ref{fig:ex2RadauOrds2357}, we can observe that the accuracy of the mean absolute errors for the 3rd and 13th-order Radau IIA methods corresponds to the blue Y-axis, while the accuracy of the average absolute errors for the 5th and 9th-order Radau methods corresponds to the red Y-axis. For the differential function variables $y_1(t)$, $y_3(t)$, and $y_4(t)$, the 3rd-order Radau IIA method has significantly higher average absolute errors than the 5th, 9th, and 13th-order methods. For the differential function variable $y_2(t)$, the 13th-order Radau IIA method exhibits significantly higher average absolute errors than the 3rd, 5th, and 9th-order methods. For the algebraic variable $y_5(t)$, the 3rd and 13th-order Radau IIA methods have significantly higher average absolute errors compared to the 5th and 9th-order methods. 

Additionally, we further observe that for all differential function variables from red Y-axis, the 9th-order Radau IIA method's average absolute error overall trends similarly to the 5th-order method. For the algebraic variable $y_5(t)$, the 9th-order Radau IIA method exhibits significantly higher average absolute errors in the later time regions compared to the 5th-order method. In other words, a PINN based on the 5th-order Radau IIA method achieves the highest precision in terms of average absolute errors.

The absolute error results obtained using the 5th-order method are shown in Figure \ref{fig:ex2RadauOrds5}. The accuracy of the absolute errors for $y_1(t)$, $y_2(t)$, $y_3(t)$, and $y_4(t)$ corresponds to the blue Y-axis, while the accuracy of the absolute errors for  $y_5(t)$ corresponds to the red Y-axis. From the figure, we can see that the lowest precision of absolute errors for all four differential variables is maintained at $10^{-7}$, while the lowest precision of absolute errors for the algebraic variable is kept at $10^{-5}$. The experimental results suggest that the neural network's predicted solutions for the pendulum's DAEs system can also achieve high precision.

% \section{Discussion}
% %%\label{}
% \lipsum[4]

\section{Summary and Conclusions}
%%\label{}
%\lipsum[1-4]
DAE systems are widely employed in various domains, including fluid dynamics, multibody dynamics, and control theory.
In practical physical modeling, most DAE models are either low-index DAEs or high-index Hessenberg-type DAEs. Classical implicit numerical methods are suitable for a certain class of high-index DAEs, but they often lead to varying order reduction of numerical accuracy.~Recently, a novel neural network method, DAE-PINN, has been developed for solving low-index DAEs. However, it cannot directly handle high-index systems. Therefore, this paper proposes a PINN-based approach using the Radau method to solve high-index DAEs systems. This method combines the strengths of the Radau IIA method with a neural network structure based on attention mechanisms and employs a time-domain decomposition strategy to enhance both efficiency and accuracy in solving these systems.

In this paper, two high-index systems, namely Hessenberg-type DAEs and pendulum model DAEs, are studied as examples. The research takes into account the influence of different orders in the Radau IIA methods and the activation functions on the accuracy of neural network solutions. Generally, employing higher-order Radau IIA methods enhances the neural network's generalization capability. However, through comparative experiments with two examples, it is found that PINN based on the 5th-order Radau IIA method provide the highest accuracy in solving the systems. This conclusion is consistent with the notion that Radau-5 is a high-precision numerical method \upcite{Ascher1998}.~Further experimental results indicate that in high-index systems, the absolute errors for all differential variables maintain a minimum precision of $10^{-6}$, while the absolute errors for algebraic variables maintain a minimum precision of $10^{-5}$. This method's numerical accuracy surpasses the corresponding results in the literature \upcite{Yang2019} and, to some extent, surpasses the accuracy achieved by the DAE-PINN method \upcite{Moya2023}. This demonstrates that our method can directly and accurately solve high-index DAEs systems, showcasing strong generalization capabilities and offering a viable approach for high-precision solutions to even higher-index DAEs or challenging systems of partial differential algebraic equations. Furthermore, we have maintained the depth and width of the neural networks as in DAE-PINN \upcite{Moya2023} and have not delved into a detailed study of their impact on the accuracy of our method, which we will need to investigate in our future work.

\section*{Acknowledgements}
Project supported by the National Natural Science Foundation of China (Grant No.~12201144), the the GuangDong Basic and Applied Basic Research Foundation of China (Grant No.~2020A1515110554), the Science and Technology Foundation of Guizhou Province (Grant No.~QKHJC-ZK[2021]YB015) of China, and Chongqing Talents Plan Youth Top-notch Project of China (Grant No.~2021000263).

%% The Appendices part is started with the command \appendix;
%% appendix sections are then done as normal sections
\appendix

%\section{Appendix title 1}
%% \label{}

%\section{Appendix title 2}
%% \label{}

%% If you have bibdatabase file and want bibtex to generate the
%% bibitems, please use
%%
%\bibliographystyle{elsarticle-harv} 
%\bibliography{example}

%% else use the following coding to input the bibitems directly in the
%% TeX file.

%%\begin{thebibliography}{00}

%% \bibitem[Author(year)]{label}
%% For example:

%% \bibitem[Aladro et al.(2015)]{Aladro15} Aladro, R., Martín, S., Riquelme, D., et al. 2015, \aas, 579, A101

%%\end{thebibliography}
\bibliographystyle{unsrt}%unsrt means cite by the appearence sequence in the paper, palin means cite by name spell 
\bibliography{main}%cite lib name
\end{document}